# A characterization of mutual absolute continuity of probability measures on a filtered space

Matthias G. Mayer

matthias.georg.mayer@gmail.com

27th November 2024

**Abstract**

We give a new characterization for mutual absolute continuity of probability measures on a filtered space. For this, we introduce a martingale limit $M$ that measures the similarity between the tails of the probability measures restricted to the filtration. The measures are mutually absolutely continuous if and only if $M = 1$ holds almost surely for both measures. In this case, the square roots of the Radon-Nikodym derivatives on the filtration converge in $L^2$. Finally, we apply the result to families of random variables and stochastic processes.

**MSC Subject Classification:** 60G; 60B10

## 1 Introduction

The purpose of this paper is to investigate the mutual absolute continuity of probability measures on a filtered space $\bigl(\Omega, \mathfrak{F}, (\mathfrak{F}_n)_{n\in\mathbb{N}}\bigr)$ and apply the characterization to families of random elements and stochastic processes. We give a review of basic definitions and known theorems in Section 2 that will also be used in the rest of this introduction. In Section 3, we give an overview of known characterizations of absolute continuity to derive characterizations of mutual absolute continuity. Let $\mathbb{P}$ and $\mathbb{P}'$ be probability measures on a filtered space $\bigl(\Omega, \mathfrak{F}, (\mathfrak{F}_n)_{n\in\mathbb{N}}\bigr)$, where $\mathfrak{F}_1 = \{\emptyset, \Omega\}$ and $\mathfrak{F} = \sigma\bigl(\bigcup_{n\in\mathbb{N}} \mathfrak{F}\bigr)$. Let $\mathbb{E}$ and $\mathbb{E}'$ denote their expectations respectively. A necessary condition for $\mathbb{P} \sim \mathbb{P}'$ is that $\mathbb{P}|_{\mathfrak{F}_n} \sim \mathbb{P}'|_{\mathfrak{F}_n}$. Assuming this condition, for all $n \in \mathbb{N}$, there is a $\mathfrak{F}_n$-measurable positive probability density $\Phi_n$, s.t. $\mathbb{P}'|\mathfrak{F}_n = \Phi_n \cdot \mathbb{P}|_{\mathfrak{F}_n}$. Now for $\varphi_n = \frac{\Phi_{n+1}}{\Phi_n}$, we see that $\varphi_{1,..,n} := \prod_{k=1}^n \varphi_k = \Phi_{n+1}$, $\mathbb{E}(\varphi_n|\mathfrak{F}_n) = 1$, and $\varphi_n$ is $\mathfrak{F}_{n+1}$-measurable.

We proceed in 4 steps.

1. For $n, k \in \mathbb{N}$, we define $M_{n,k} = \mathbb{E}\bigl(\sqrt{\varphi_{n,..,k}} \mid \mathfrak{F}_n\bigr)$. Morally, $M_{n,k}$ can be viewed as a $L^2$-valued inner product between the integrated kernels from $\mathfrak{F}_n$ to $\mathfrak{F}_k$ of $\mathbb{P}$ and $\mathbb{P}'$. We show that $M_n = \lim_k M_{n,k}$ and $M = \lim_n M_n$ exist $(\mathbb{P} + \mathbb{P}')$-almost everywhere. Morally, $M$ can be viewed as the similarity of $\mathbb{P}$ and $\mathbb{P}'$ at infinity. In other words, $M$ measures the similarity between $\mathbb{P}$ and $\mathbb{P}'$ that does not depend on $\mathbb{P}|_{\mathfrak{F}_n}$ nor $\mathbb{P}'|_{\mathfrak{F}_n}$ for any $n \in \mathbb{N}$. Similarly, we define $N_n = \lim_k \mathbb{E}\bigl(\sqrt{\varphi_{1,..,k}} \mid \mathfrak{F}_n\bigr) = \sqrt{\varphi_{1,..,n-1}} M_n$ and show that $N = \lim_n N_n$ exists $\mathbb{P}$-almost everywhere. Morally, $N^2$ is the limit of $\varphi_{1,...,n}$ where the effects of the tail of $(\varphi_n)_{n\in\mathbb{N}}$ are included. Finally, we introduce the same definitions for $\mathbb{P}'$, which we denote by $M'_{n,k}, M'_n, M', N'_n$ and $N'$.

2. We study the interaction between $N$ and $M$ and their properties. In particular, we have $\mathbb{P}(N > 0, M = 0) = 0$ and conclude $M(\omega) = 1$ for $\mathbb{P}$-a.e. $\omega \in \{N > 0\}$.
3. We show that $\mathbb{P}|_{\{N=0 \vee N'=0\}} \perp \mathbb{P}'|_{\{N=0 \vee N'=0\}}$.
4. We show that if $M = 1$ $(\mathbb{P} + \mathbb{P}')$-almost surely, then $\sqrt{\varphi_{1,\dots,n}}$ converges in $L^2(\mathbb{P})$. Therefore, we conclude $\mathbb{P} \sim \mathbb{P}'$.

Finally, we conclude with the full statement of the characterization and apply it to families of random variables.

## 2 Background

This section is intended to disambiguate notation. We don't provide proofs of the commonly known theorems. They can be found, for example, in [1]. Furthermore, we will apply the theorems without explicit reference to their appearance here. Let $\mathfrak{P}(\Omega)$ denote the powerset of $\Omega$. and for $S \subseteq \mathfrak{P}(\Omega)$, let $\sigma(S)$ denote the $\sigma$-algebra generated by $S$.

**Definition 2.1** (measurable space): Let $\Omega$ be a set and $\mathcal{A} \subseteq \mathfrak{P}(\Omega)$ a $\sigma$-algebra. We call $(\Omega, \mathcal{A})$ a measurable space.

**Definition 2.2** (random element): Let $(\Omega, \mathcal{A}), (\Omega', \mathcal{A}')$ be measurable spaces. We call a function $X : \Omega \to \Omega'$ that is $\mathcal{A}$-$\mathcal{A}'$-measurable a random element from $(\Omega, \mathcal{A})$ to $(\Omega', \mathcal{A}')$. It is convenient to introduce $X$ as a random element without referring explicitly to $(\Omega, \mathcal{A})$ or $(\Omega', \mathcal{A}')$, when $(\Omega, \mathcal{A})$ is understood from context. $(\Omega', \mathcal{A}')$ is denoted by $(\mathrm{Val}(X), \mathcal{V}(X))$.

**Definition 2.3** (Dynkin system): A Dynkin system on $\Omega$ is a set $\mathcal{D} \subseteq \mathfrak{P}(\Omega)$, s.t. $\Omega \in \mathcal{D}$, and $A, B \in \mathcal{D} : A \subseteq B \Rightarrow B \setminus A \in \mathcal{D}$, and for $(A_n)_{n \in \mathbb{N}} \in \mathcal{D}^{\mathbb{N}}$ pairwise disjoint, $\bigcup_{n \in \mathbb{N}} A_n \in \mathcal{D}$.

**Definition 2.4** ($\cap$-stable system): A $\cap$-stable system (or $\pi$-system), is a set $\mathcal{B} \subseteq \mathfrak{P}(\Omega)$, s.t. $A, B \in \mathcal{B} \Rightarrow A \cap B \in \mathcal{B}$.

**Theorem 2.5**: If $S$ is a $\cap$-stable system and $\mathcal{D}$ a Dynkin system, then $S \subseteq \mathcal{D} \Rightarrow \sigma(S) \subseteq \mathcal{D}$.

**Theorem 2.6** (generator approximation): Let $\mu$ be a finite measure on $(\Omega, \mathcal{A})$. Let $S$ be a generators of $\mathcal{A}$. Suppose that $S$ is $\cap$-stable and $A \in S \Rightarrow A^c \in S$. Then for any $A \in \mathcal{A}$ and $\varepsilon > 0$, there is $B \in S$, s.t. $\mu(A \triangle B) < \varepsilon$.

**Definition 2.7** (product measurable space): Given a family of measurable spaces $(\Omega_i, \mathcal{A}_i)_{i \in I}$, the product measurable space is defined by $\bigotimes_{i \in I}(\Omega_i, \mathcal{A}_i) := (\Omega, \mathcal{A})$, where $\Omega = \bigtimes_{i \in I}$ and $\mathcal{A}$ is the $\sigma$-algebra on $\Omega$ that renders the projections $\pi_i : \Omega \to \Omega_i$ measurable.

**Lemma 2.8**: Given a product space $(\Omega, \mathcal{A}) = \bigotimes_{i \in I}(\Omega_i, \mathcal{A}_i)$, for all $A \in \mathcal{A}$, there is a countable set $I_0 \subseteq I$, s.t. $A \in \sigma(\pi_i : i \in I_0)$.

**Definition 2.9**: Given a family of random elements $(X_i)_{i \in I}$ defined on a common measurable space $(\Omega, \mathcal{A})$, we associate to $(X_i)_{i \in I}$ the random element defined by $(\Omega, \mathcal{A}) \to \bigotimes_{i \in I}(\mathrm{Val}(X_i), \mathcal{V}(X_i)); \omega \mapsto (X_i(\omega))_{i \in I}$.

Let $(\Omega, \mathcal{A})$ be a measurable space and $\mathbb{P}$ and $\mathbb{P}'$ be probability distributions on this space. Let $\mathbb{E}$ and $\mathbb{E}'$ denote their expectations.

**Definition 2.10** (density): A measurable function $\varphi : \Omega \to [0, \infty]$ is called density. It is called probability density w.r.t. $\mathbb{P}$, if $\mathbb{E}(\varphi) = 1$. We define the measure $\varphi \cdot \mathbb{P}$ by $(\varphi \cdot \mathbb{P})(A) := \mathbb{E}(1_A \varphi)$.

**Definition 2.11** (orthogonality): We call $\mathbb{P}$ and $\mathbb{P}'$ orthogonal and denote this by $\mathbb{P} \perp \mathbb{P}'$, if there is a $C \in \mathcal{A}$, s.t. $\mathbb{P}(C) = 0$ and $\mathbb{P}'(\Omega \setminus C) = 0$.

**Definition 2.12** (absolute continuity): We call $\mathbb{P}'$ absolutely continuous w.r.t. $\mathbb{P}$ and denote this by $\mathbb{P}' \ll \mathbb{P}$, if $\forall A \in \mathcal{A} : \mathbb{P}(A) = 0 \Rightarrow \mathbb{P}'(A) = 0$.

**Definition 2.13** (mutual absolute continuity): We write $\mathbb{P} \sim \mathbb{P}'$, if $\mathbb{P} \ll \mathbb{P}'$ and $\mathbb{P}' \ll \mathbb{P}$.

**Theorem 2.14** (Radon-Nikodym derivative): If $\mathbb{P}' \ll \mathbb{P}$, there exists a density $\varphi$ w.r.t. $\mathbb{P}$, s.t. $\varphi \cdot \mathbb{P} = \mathbb{P}'$. We write $\varphi =: \frac{d\mathbb{P}'}{d\mathbb{P}}$. If $\mathbb{P}' \sim \mathbb{P}$, then $\varphi$ can be chosen to be $\mathbb{P}$-a.s. positive.

**Definition 2.15** (conditional expectation): Let $X \in L^1(\mathbb{P})$ and $\mathcal{C}$ a sub-$\sigma$-algebra of $\mathcal{A}$. There exists an up to $\mathbb{P}$-nullsets unique, $\mathcal{C}$-measurable map, $\mathbb{E}(X|Z)$ that fulfills $\mathbb{E}(1_C X) = \mathbb{E}(1_C \mathbb{E}(X|Z))$ for all $C \in \mathcal{C}$.

**Definition 2.16** (filtration): A sequence of $\sigma$-algebras $(\mathfrak{F}_n)_{n \in \mathbb{N}}$ on $\Omega$ is called a filtration if $n < m \Rightarrow \mathfrak{F}_n \subseteq \mathfrak{F}_m$.

**Definition 2.17** ((sub/super)martingale): A stochastic process $(X_n)_{n \in \mathbb{N}}$ is called a martingale w.r.t. the filtration $(\mathfrak{F}_n)_{n \in \mathbb{N}}$ and $\mathbb{P}$, if for all $n \in \mathbb{N}$, $X_n$ is integrable and $\mathfrak{F}_n$-measurable and $X_n = \mathbb{E}(X_{n+1} \mid \mathfrak{F}_n)$. It is called sub/supermartingale if equality is replaced by $\leq$ or $\geq$ respectively. For $p > 1$, $X_n$ is called a (uniform) $L^p$-(sub/super)martingale, if $X_n \in L^p(\mathbb{P})$ ($\sup_{n \in \mathbb{N}} \|X_n\|_{L^p(\mathbb{P})} < \infty$).

The following are variants of the Doob's martingale convergence theorems.

**Theorem 2.18**: If $(X_n)_{n \in \mathbb{N}}$ is a supermartingale w.r.t. $\mathbb{P}$ that is bounded from below, then the point-wise limit $X = \lim_n X_n$ exists (finitely) $\mathbb{P}$-almost everywhere. The same is true if $X_n$ is a submartingale bounded from above. If $(X_n)_{n \in \mathbb{N}}$ is a uniform $L^p$-martingale for $p > 1$ Then $X_n$ converges in $L^p$ and $\mathbb{P}$-a.s.

## 3 Related Work

Let $\mathbb{P}$ and $\mathbb{P}'$ be defined on a filtered space $\left(\Omega, \mathfrak{F}, (\mathfrak{F}_n)_{n \in \mathbb{N}}\right)$, where $\mathfrak{F}_1 = \{\emptyset, \Omega\}$ and $\mathfrak{F} = \sigma\left(\bigcup_{n \in \mathbb{N}} \mathfrak{F}_n\right)$, s.t. $\mathbb{P}|_{\mathfrak{F}_n} \sim \mathbb{P}'|_{\mathfrak{F}_n}$. The question of whether $\mathbb{P}' \ll \mathbb{P}$ holds in this case has been answered in two ways.

Let $\Phi_n := \frac{d\mathbb{P}'|_{\mathfrak{F}_n}}{d\mathbb{P}|_{\mathfrak{F}_n}}$. Kabanov et al. [2], Lemma 6, and [3], show that $\Phi = \lim_n \Phi_n$ exists $\mathbb{P}'$-a.s. Furthermore, $\mathbb{P}$ is orthogonal to $\mathbb{P}'$ on the set $\{\Phi = \infty\}$ and absolutely continuous on $\{\Phi < \infty\}$. More precisely, $\mathbb{P}'(A) = \int_A \Phi \, d\mathbb{P} + \mathbb{P}'(A \cap \{\Phi = \infty\})$. Therefore, $\mathbb{P} \ll \mathbb{P}' \Leftrightarrow \Phi < \infty$ $\mathbb{P}'$-a.s.. Reversing the roles of $\mathbb{P}$ and $\mathbb{P}'$, we obtain that $\Phi' = \lim_n \Phi_n^{-1}$ exists $\mathbb{P}$-a.s. Therefore $\Phi = \lim_n \Phi_n$ exists $(\mathbb{P} + \mathbb{P}')$-a.e. and $\mathbb{P} \sim \mathbb{P}' \Leftrightarrow \Phi < \infty$ $\mathbb{P}'$-a.s. and $\Phi > 0$ $\mathbb{P}$-a.s.

We can immediately collapse the conditions on $\mathbb{P}$ and $\mathbb{P}'$.

**Corollary 3.1**: Let $\mathbb{P}$ and $\mathbb{P}'$ be probability distributions on a filtered space $(\Omega, \mathfrak{F}, \mathfrak{F}_{n \in \mathbb{N}})$, s.t. $\mathbb{P}|_{\mathfrak{F}_n} \sim \mathbb{P}'|_{\mathfrak{F}_n}$. Let $\Phi_n = \frac{d\mathbb{P}'|_{\mathfrak{F}_n}}{d\mathbb{P}|_{\mathfrak{F}_n}}$ and $\Phi = \lim_n \Phi_n$ $(\mathbb{P} + \mathbb{P}')$-a.e. Then $\mathbb{P} \sim \mathbb{P}' \Leftrightarrow 0 < \Phi < \infty$ $(\mathbb{P} + \mathbb{P}')$-a.e.

*Proof*: Since $\mathbb{P} \sim \mathbb{P}' \Leftrightarrow \Phi < \infty$ $\mathbb{P}'$-a.s. and $\Phi > 0$ $\mathbb{P}$-a.s. '$\Leftarrow$' is a direct implication of this observation. '$\Rightarrow$': Let $\mathbb{P} \sim \mathbb{P}'$, then $\Phi < \infty$ $\mathbb{P}'$-a.s. and since $\mathbb{P} \sim \mathbb{P}'$, also $\Phi < \infty$ $\mathbb{P}$-a.s. Similarly, $\Phi > 0$ $(\mathbb{P} + \mathbb{P}')$-a.s. □

Kabunov et al. uses this criterion for absolute continuity to introduce a similar condition to Kakutani's result on infinite product probability measures [4]. Kakutani introduces a distance between two probability measures by $d(\mathbb{P}', \mathbb{P}) = \int \sqrt{\frac{d\mathbb{P}'}{d\mathbb{P}}} d\mathbb{P}$. Then two infinite product probability distributions $\mathbb{P} = \bigtimes_{n \in \mathbb{N}} \mathbb{P}_n$ and $\mathbb{P}' = \bigtimes_{n \in \mathbb{N}} \mathbb{P}'_n$ are mutually absolutely continuous if and only if $\prod_{n \in \mathbb{N}} d(\mathbb{P}_n, \mathbb{P}'_n) > 0$. This condition is equivalent to $\sum_{n \in \mathbb{N}} |\ln d(\mathbb{P}_n, \mathbb{P}'_n)| < \infty$.

Kabunov et al. produces a similar result. Let $\varphi_n = \Phi_{n+1}\Phi_n^{-1}$. Then $\mathbb{P} \ll \mathbb{P}'$ if and only if $\sum_{n \in \mathbb{N}} |\ln \mathbb{E}(\sqrt{\varphi_n} \mid \mathfrak{F}_n)| < \infty$ $\mathbb{P}'$-a.s. Therefore we have $\mathbb{P} \sim \mathbb{P}'$, if and only if $\sum_{n \in \mathbb{N}} |\ln \mathbb{E}(\sqrt{\varphi_n} \mid \mathfrak{F}_n)| < \infty$ $\mathbb{P}'$-a.s. and $\sum_{n \in \mathbb{N}} |\ln \mathbb{E}'(\sqrt{\varphi_n^{-1}} \mid \mathfrak{F}_n)| < \infty$ $\mathbb{P}$-a.s.

Like above, we can collapse the requirements.

**Corollary 3.2**: Let $\mathbb{P}$ and $\mathbb{P}'$ be probability distributions on a filtered space $(\Omega, \mathfrak{F}, \mathfrak{F}_{n \in \mathbb{N}})$, s.t. $\mathbb{P}|_{\mathfrak{F}_n} \sim \mathbb{P}'|_{\mathfrak{F}_n}$. Let $\Phi_n = \frac{d\mathbb{P}'|_{\mathfrak{F}_n}}{d\mathbb{P}|_{\mathfrak{F}_n}}$ and $\varphi_n = \Phi_{n+1}\Phi_n^{-1}$. Then $\mathbb{P} \sim \mathbb{P}' \Leftrightarrow \sum_{n \in \mathbb{N}} |\ln \mathbb{E}(\sqrt{\varphi_n} \mid \mathfrak{F}_n)| < \infty$ $(\mathbb{P} + \mathbb{P}')$-a.e.

*Proof*: Trivial from the discussion above. □

The result in this paper was discovered independently from Kabunov et al. and is different in some aspects. A key difference is that this work stays closer to Kakutani's proof method and establishes that $\sqrt{\Phi_n}$ converges to $\sqrt{\frac{d\mathbb{P}'}{d\mathbb{P}}}$ in $L^2(\mathbb{P})$ in case of mutual absolute continuity.

While Kabunov et al's criterion is based on a predictable convergence, i.e. the limit takes terms that are $\mathfrak{F}_n$-measurable, $M$ is not a limit of a predictable sequence. Therefore the criterion $M = 1$ $(\mathbb{P} + \mathbb{P}')$-a.e. is of a more theoretical nature.

# 4 A characterization of mutual absolute continuity

Let $(\Omega, \mathfrak{F}, (\mathfrak{F}_n)_{n \in \mathbb{N}})$ be a filtered space, s.t. $\mathfrak{F}_1 = \{\emptyset, \Omega\}$ and $\mathfrak{F} = \sigma(\bigcup_{n \in \mathbb{N}} \mathfrak{F}_n)$. Let $\mathbb{P}$ and $\mathbb{P}'$ be probability measures on $(\Omega, \mathfrak{F})$, s.t. $\mathbb{P}|_{\mathfrak{F}_n} \sim \mathbb{P}'|_{\mathfrak{F}_n}$ for all $n \in \mathbb{N}$. Our goal is to investigate when $\mathbb{P} \sim \mathbb{P}'$ holds.

Let $\Phi_n = \frac{d\mathbb{P}'|_{\mathfrak{F}_n}}{d\mathbb{P}|_{\mathfrak{F}_n}}$ and Let $\Phi'_n = \frac{d\mathbb{P}|_{\mathfrak{F}_n}}{d\mathbb{P}'|_{\mathfrak{F}_n}}$. Note that $\Phi_n^{-1} = \Phi'_n$. For $n \in \mathbb{N}$, let $\varphi_n = \frac{\Phi_{n+1}}{\Phi_n}$ and $\varphi'_n = \frac{\Phi'_{n+1}}{\Phi'_n}$. Since for any $A \in \mathfrak{F}_n$, $\mathbb{E}(1_A \Phi_{n+1}) = \mathbb{E}'(1_A) = \mathbb{E}(1_A \Phi_n)$, we have $\mathbb{E}(\Phi_{n+1} \mid \mathfrak{F}_n) = \Phi_n$ by the defining property of conditional expectation. Therefore $\mathbb{E}(\varphi_n \mid \mathfrak{F}_n) = 1$. Note that $\varphi_n^{-1} = \varphi'_n$.

During the next section, we will introduce the following definitions. There will be some work showing that these are well defined. For $n, k \in \mathbb{N}$, let $M_{n,k} := \mathbb{E}(\sqrt{\varphi_{n,..,k}} \mid \mathfrak{F}_n)$. Note that $M_{n,k}$ is $\mathfrak{F}_k$-measurable, so that it is defined up $(\mathbb{P} + \mathbb{P}')$-almost everywhere. Similarly, we define $M'_{n,k} = \mathbb{E}'(\sqrt{\varphi'_{n,..,k}} \mid \mathfrak{F}_n)$ and will see that $M_{n,k} = M'_{n,k}$.

Our goal in the next section is to show that $M_n := \lim_k M_{n,k}$, and $M := \lim_n M_n$ exists $(\mathbb{P} + \mathbb{P}')$-almost everywhere. Then we define $N_n := \sqrt{\varphi_{1,..,n-1}} M_n$ and show that $N := \lim_n N_n$ exists $\mathbb{P}$-a.s. Similarly, we define $N'_n := \sqrt{\varphi'_{1,..,n-1}} M'_n$ and $N' := \lim_n N'_n$.

## 4.1 Convergence

**Definition 4.1.1**: For $n, k \in \mathbb{N}$, let $M_{n,k} := \mathbb{E}(\sqrt{\varphi_{n,...,k}} \mid \mathfrak{F}_n)$ and $M'_{n,k} := \mathbb{E}'(\sqrt{\varphi'_{n,...,k}} \mid \mathfrak{F}_n)$. Note that by $\mathfrak{F}_n$ measurability and $\mathbb{P}|_{\mathfrak{F}_n} \sim \mathbb{P}'|_{\mathfrak{F}_n}$, we have that $M_{n,k}$ and $M'_{n,k}$ are defined up to $(\mathbb{P} + \mathbb{P}')$-nullsets.

**Lemma 4.1.2**: $M_{n,k} = M'_{n,k}$ $(\mathbb{P} + \mathbb{P}')$-almost everywhere.

*Proof*: We use the defining property of conditional expectation. Clearly, $M'_{n,k}$ is $\mathfrak{F}_n$ measurable. Furthermore, for $A \in \mathfrak{F}_n$, we have $\int_A \mathbb{E}'\left(\sqrt{\varphi'_{n,..,k}} \mid \mathfrak{F}_n\right) d\mathbb{P} = \int_A \varphi'_{1,..,n-1} \mathbb{E}'\left(\sqrt{\varphi'_{n,..,k}} \mid \mathfrak{F}_n\right) d\mathbb{P}' = \int_A \varphi'_{1,...,n-1} \sqrt{\varphi'_{n,..,k}} d\mathbb{P}' = \int_A \sqrt{\varphi_{n,..,k}} d\mathbb{P}$. $\square$

**Lemma 4.1.3**: For $n, k \in \mathbb{N}$, we have $1 \geq M_{n,k} \geq M_{n,k+1} \geq 0$ $(\mathbb{P} + \mathbb{P}')$-almost everywhere.

*Proof*: Since $\mathbb{P}|_{\mathfrak{F}_{k+1}} \sim \mathbb{P}'|_{\mathfrak{F}_{k+1}}$, it suffices to show the inequality $\mathbb{P}$-almost surely. By Jensen's inequality, $M_n^k = \mathbb{E}(\sqrt{\varphi_{n,...,k}} \mid \mathfrak{F}_n) \leq \sqrt{\mathbb{E}(\varphi_{n,...,k} \mid \mathfrak{F}_n)} = 1$. Furthermore, $M_{n,k+1} = \mathbb{E}\left(\sqrt{\varphi_{n,...,k+1}} \mid \mathfrak{F}_n\right) = \mathbb{E}(\sqrt{\varphi_{n,...,k}} \cdot M_{k+1,k+1} \mid \mathfrak{F}_n) \leq M_{n,k}$. $\square$

**Lemma 4.1.4**: $M_{n,k}$ converges $(\mathbb{P} + \mathbb{P}')$-almost everywhere for $k \to \infty$.

*Proof*: Follows directly from Lemma 4.1.3. $\square$

**Definition 4.1.5**: For $n \in \mathbb{N}$, let $M_n := \lim_k M_{n,k}$.

**Lemma 4.1.6**: $M_n^2$ is a submartingale w.r.t. the filtration $(\mathfrak{F}_n)_{n \in \mathbb{N}}$ and $\mathbb{P}$. The same holds w.r.t. $\mathbb{P}'$.

*Proof*: We first show the case w.r.t. $\mathbb{E}$. Note that by dominated convergence, $\mathbb{E}(M_n^2) = \mathbb{E}\left(\lim_k \mathbb{E}\left(\sqrt{\varphi_{n,...,k}} \mid \mathfrak{F}_n\right)^2\right) \leq \lim_k \mathbb{E}(\varphi_{n,...,k})) = 1$. Furthermore, $M_n^2$ is $\mathfrak{F}_n$-measurable. Next, by Hölder's inequality for conditional expectation, $M_n = \lim_k \mathbb{E}(\sqrt{\varphi_{n,...k}} \mid \mathfrak{F}_n) = \lim_k \mathbb{E}(\sqrt{\varphi_n} \cdot M_{n+1,k} \mid \mathfrak{F}_n) \leq \lim_k (\mathbb{E}(\varphi_n \mid \mathfrak{F}_n) \cdot \mathbb{E}(M_{n+1,k}^2 \mid \mathfrak{F}_n)))^{\frac{1}{2}} = \lim_k \mathbb{E}(M_{n+1,k}^2 \mid \mathfrak{F}))^{\frac{1}{2}}$. By dominated convergence, the right hand side equals $\mathbb{E}(M_{n+1}^2 \mid \mathfrak{F}_n)^{\frac{1}{2}}$. Therefore, $M_n^2 \leq \mathbb{E}(M_{n+1}^2 \mid \mathfrak{F}_n)$.

By the symmetry between $\mathbb{P}$ and $\mathbb{P}'$, we see that $M'^2_n$ is a submartingale w.r.t $\mathbb{E}'$. Since $M_n^2 = M'^2_n$ $(\mathbb{P} + \mathbb{P}')$-a.e. the case w.r.t. $\mathbb{E}'$ follows. $\square$

**Lemma 4.1.7**: $M_n$ converges $(\mathbb{P} + \mathbb{P}')$-almost everywhere

*Proof*: $M_n^2$ is a bounded submartingale w.r.t. $\mathbb{P}$ and $\mathbb{P}'$, therefore it converges $\mathbb{P}$- and $\mathbb{P}'$-almost surely. Therefore $M_n$ converges $(\mathbb{P} + \mathbb{P}')$-almost surely. $\square$

**Definition 4.1.8**: For $n \in \mathbb{N}$, let $N_n := \sqrt{\varphi_{1,...,n-1}} \cdot M_n = \lim_k \mathbb{E}\left(\sqrt{\varphi_{1,...,k}} \mid \mathfrak{F}_n\right)$, and $N'_n = \sqrt{\varphi'_{1,...,n-1}} \cdot M'_n$.

**Lemma 4.1.9**: $N_n$ is a uniform $L^2$-martingale w.r.t. the filtration $(\mathfrak{F}_n)_{n \in \mathbb{N}}$ and $\mathbb{P}$.

*Proof*: First, $\mathbb{E}(N_n^2) = \mathbb{E}\left(\varphi_{1,...,n-1} \lim_k \mathbb{E}\left(\sqrt{\varphi_{n,...,k}} \mid \mathfrak{F}_n\right)^2\right) \leq \mathbb{E}(\varphi_{1,...,n-1}) = 1$. Clearly, $N_n$ is $\mathfrak{F}_n$-measurable. By dominated convergence with dominant $\sqrt{\varphi_{1,...,n}}$ and the tower property, we have $\mathbb{E}(N_{n+1}|\mathfrak{F}_n) = \mathbb{E}\left(\lim_k \mathbb{E}\left(\sqrt{\varphi_{1,...,k}} \mid \mathfrak{F}_{n+1}\right) \mid \mathfrak{F}_n\right) = \lim_k \mathbb{E}\left(\mathbb{E}\left(\sqrt{\varphi_{1,...,k}} \mid \mathfrak{F}_{n+1}\right) \mid \mathfrak{F}_n\right) = \lim_k \mathbb{E}\left(\sqrt{\varphi_{1,...,k}} \mid \mathfrak{F}_n\right) = N_n$. $\square$

**Corollary 4.1.10**: $N_n$ converges $\mathbb{P}$-almost surely and in $L^2(\mathbb{P})$. Similarly, $N'_n$ converges $\mathbb{P}'$-almost surely and in $L^2(\mathbb{P}')$.

*Proof*: Follows immediately from Lemma 4.1.9 and martingale convergence theorems. $\square$

**Definition 4.1.11**: Let $N = \lim_n N_n$ and $N' = \lim_n N'_n$, where we take the limit $\mathbb{P}$ and $\mathbb{P}'$ almost surely respectively.

## 4.2 The relation between $N$ and $M$

**Lemma 4.2.1**: $\mathbb{P}(N > 0, M = 0) = 0$.

*Proof*: First, note that for $\omega \in \{N > 0, M = 0\} =: B$, we have that $\sqrt{\varphi_{1,...,n-1}}(\omega) \cdot M_n(\omega) \to N(\omega) > 0$. Since $M_n(\omega) \to M(\omega) = 0$, $\sqrt{\varphi_{1,...k}} \to \infty$. Now by Fatou's lemma,

$$\int_B \liminf_k \sqrt{\varphi_{n,\ldots,k}} \, d\mathbb{P} \leq \liminf_k \int_B \sqrt{\varphi_{n,\ldots,k}} \, d\mathbb{P}$$
$$\leq \liminf_k \int \sqrt{\varphi_{n,\ldots,k}} \, d\mathbb{P} \tag{1}$$
$$= \liminf_k \int M_{n,k} \, d\mathbb{P}$$
$$\leq 1$$

Therefore $\int \infty \cdot 1_B \, d\mathbb{P} < \infty$, so $\mathbb{P}(B) = 0$. $\square$

**Lemma 4.2.2**: $\varphi_{1,\ldots,n}(\omega)$ converges to a real number for $n \to \infty$ and $\mathbb{P}$-a.e. $\omega \in \{M > 0\}$.

*Proof*: We have $N_n \to N$ and $M_n \to M$ $\mathbb{P}$-a.s. Furthermore, note that $N_n = \sqrt{\varphi_{1,\ldots,n-1}} \cdot M_n$. Therefore, $\sqrt{\varphi_{1,\ldots,n-1}}(\omega) = \frac{N_n(\omega)}{M_n(\omega)} \to \frac{N(\omega)}{M(\omega)} \in \mathbb{R}$ for a.e. $\omega \in \{M > 0\}$. $\square$

**Lemma 4.2.3**: $\varphi_{1,\ldots,n}(\omega)$ converges to a positive real number for $n \to \infty$ and $\mathbb{P}$-a.e. $\omega \in \{N > 0\}$.

*Proof*: We have $N_n \to N$ and $M_n \to M$ $\mathbb{P}$-a.s. By Lemma 4.2.1, for a.e. $\omega \in \{N > 0\}$ we have $M(\omega) > 0$. Furthermore, note that $N_n = \sqrt{\varphi_{1,\ldots,n-1}} \cdot M_n$. Therefore, $\sqrt{\varphi_{1,\ldots,n-1}}(\omega) = \frac{N_n(\omega)}{M_n(\omega)} \to \frac{N(\omega)}{M(\omega)} \in \mathbb{R}_{>0}$ for a.e. $\omega \in \{N > 0\}$. $\square$

**Lemma 4.2.4**: $\lim_n \lim_k \varphi_{n,\ldots k}(\omega) \to 1$ for $\mathbb{P}$-a.e. $\omega \in \{N > 0\}$.

*Proof*: Since $\varphi_{1,\ldots,k}(\omega) = \prod_{n=1}^k \varphi_n(\omega)$ converges to a positive real number for $\mathbb{P}$-a.e. $\omega \in \{N > 0\}$, the claim follows. $\square$

**Theorem 4.2.5**: $M(\omega) = 1$ for $\mathbb{P}$-a.e. $\omega \in \{N > 0\}$ and for $\mathbb{P}'$-a.e. $\omega \in \{N' > 0\}$.

*Proof*: By symmetry, it suffices to show the claim for $C = \{N > 0\}$ and $\mathbb{P}$.

For $n \in \mathbb{N}$, let $\varphi_{n,\ldots} = \liminf_k \varphi_{n,\ldots,k}$. Let $n \in \mathbb{N}$ and $D \in \mathfrak{F}_n$. By Fatou's lemma, $\mathbb{E}\left(\sqrt{\varphi_{n,\ldots}} \mid \mathfrak{F}_n\right) \leq \lim_k \mathbb{E}\left(\sqrt{\varphi_{n,\ldots,k}} \mid \mathfrak{F}_n\right) = M_n$. Let $a \wedge b := \min(a,b)$. Then
$$\mathbb{E}\left(1_D \liminf_n \sqrt{\varphi_{n,\ldots}} \wedge 1\right) \leq \liminf_n \mathbb{E}\left(1_D \sqrt{\varphi_{n,\ldots}} \wedge 1\right)$$
$$= \liminf_n \mathbb{E}\left(1_D \mathbb{E}\left(\sqrt{\varphi_{n,\ldots}} \wedge 1 \mid \mathfrak{F}_n\right)\right) \tag{2}$$
$$= \lim_n \mathbb{E}(1_D M_n)$$
$$= \mathbb{E}(1_D M)$$

where we used dominated convergence. Since the sets $D$ that fulfill $\mathbb{E}\left(1_D \liminf_n \sqrt{\varphi_{n,\ldots}} \wedge 1\right) \leq \mathbb{E}(1_D M)$ form a Dynkin system and the inequality holds for all $D \in \bigcup_{n \in \mathbb{N}} \mathfrak{F}_n$, the inequality holds for all $D \in \mathfrak{F}$ and in particular for $C$.

By the preceding lemmas, we have $\varphi_{n,\ldots}(\omega) = 1$ for a.e. $\omega \in C$. Therefore, $\mathbb{E}(1_C) \leq \mathbb{E}(1_C M)$. Since $M \leq 1$, we have $M(\omega) = 1$ for a.e. $\omega \in C$. $\square$

## 4.3 Orthogonality

**Theorem 4.3.1**: For $C \in \{\{N = 0\}, \{N' = 0\}\}$, we have $\mathbb{P}|_C \perp \mathbb{P}'|_C$.

*Proof*: By symmetry, it suffices to show the claim for $C = \{N = 0\}$. Let $\varepsilon > 0$. Let $n \in \mathbb{N}$, s.t. $\int_C N_n \leq \int_C N + \varepsilon = \varepsilon$ ($L^1$ convergence), and $C' \in \mathfrak{F}_n$, s.t. $\mathbb{P}(C' \triangle C) < \varepsilon$ and $\mathbb{P}'(C' \triangle C) < \varepsilon$ (generator approximation theorem), and $\forall m \in \mathbb{N}: \int_{C'} N_m \, d\mathbb{P} \leq \int_C N_m \, d\mathbb{P} + \varepsilon$. (uniform integrability).

Let $k > n$, s.t. $\int_{C'} \sqrt{\varphi_{1,\ldots,n-1}} \cdot M_{n,k}\,d\mathbb{P} \leq \int_{C'} N_n\,d\mathbb{P} + \varepsilon$. ($L^1$ convergence, by dominated convergence)

Let $B = \{\varphi_{1,\ldots,k} > 1\}$. Then
$$\mathbb{P}(C' \cap B) \leq \int_{C' \cap B} \sqrt{\varphi_{1,\ldots,k}}\,d\mathbb{P} \leq \int_{C'} \sqrt{\varphi_{1,\ldots,k}}\,d\mathbb{P} = \int_{C'} \sqrt{\varphi_{1,\ldots,n-1}}M_{n,k}\,d\mathbb{P} \tag{3}$$

Similarly,
$$\mathbb{P}'(C' \setminus B) = \int_{C' \setminus B} \varphi_{1,\ldots,k}\,d\mathbb{P} \leq \int_{C' \setminus B} \sqrt{\varphi_{1,\ldots,k}}\,d\mathbb{P} \leq \int_{C'} \sqrt{\varphi_{1,\ldots,n-1}}M_{n,k}\,d\mathbb{P} \tag{4}$$

By the choices of $n$ and $k$, we have
$$\int_{C'} \sqrt{\varphi_{1,\ldots,n-1}}M_{n,k} \leq \int_{C'} N_n\,d\mathbb{P} + \varepsilon \leq \int_C N_n + 2\varepsilon \leq 3\varepsilon \tag{5}$$

Therefore, $\mathbb{P}(C' \cap B) < 3\varepsilon$ and $\mathbb{P}'(C' \setminus B) \leq 3\varepsilon$. Since $\mathbb{P}(C' \triangle C) < \varepsilon$ and $\mathbb{P}'(C' \triangle C) < \varepsilon$, we have $\mathbb{P}(C \cap B) < 4\varepsilon$ and $\mathbb{P}'(C \setminus B) \leq 4\varepsilon$. $\square$

**Corollary 4.3.2**: Suppose $M = 1$ does not hold $(\mathbb{P} + \mathbb{P}')$-a.e. Then $\mathbb{P} \not\sim \mathbb{P}'$.

*Proof*: By symmetry, we can assume that there is a set $C$, s.t. $\mathbb{P}(C) > 0$ and $\forall \omega \in C : M(\omega) < 1$. Then by Theorem 4.2.5, we have $N(\omega) = 0$ for a.e. $\omega \in C$. Therefore $\mathbb{P}(\{N = 0\}) > 0$ and by Theorem 4.3.1, $\mathbb{P}|_C \perp \mathbb{P}'|_C$. Since $\mathbb{P}(C) \neq 0$, this implies $\mathbb{P}|_C \not\sim \mathbb{P}'|_C$. $\square$

## 4.4 Mutual absolute continuity

**Lemma 4.4.1**: Let $(S, d)$ be a complete metric space. Let $(s_n)_{n \in \mathbb{N}}$ be a sequence in $S$. Suppose that $\lim_n \lim_k d(s_n, s_k) \to 0$ and $d(s_n, s_k) \leq d(s_n, s_{k+1})$ for all $n, k \in \mathbb{N}$. Then $s_n$ converges.

*Proof*: We show that $s_n$ is a Cauchy sequence. Let $\varepsilon > 0$. Let $d_\infty(s_n) := \lim_m d(s_n, s_m)$. Note that by assumption $\forall m \in \mathbb{N}, d_\infty(s_n) \geq d(s_n, s_m)$. Choose $n$, s.t. $\forall m > n : d_\infty(s_m) < \varepsilon$. Then for $n < m < k$, we have $d(s_m, s_k) \leq d_\infty(s_m) \leq \varepsilon$. $\square$

**Lemma 4.4.2**: Let $M = 1$ $(\mathbb{P} + \mathbb{P}')$-a.e. Then $\sqrt{\varphi_{1,\ldots,n}}$ converges in $L^2(\mathbb{P})$ and $(\mathbb{P} + \mathbb{P}')$-a.e.

*Proof*: By Lemma 4.2.2, we have that $\varphi_{1,\ldots,n}$ converges $(\mathbb{P} + \mathbb{P}')$-a.e. Since $L^2(\mathbb{P})$ is a Hilbert space, we have
$$\left\|\sqrt{\varphi_{1,\ldots,n-1}} - \sqrt{\varphi_{1,\ldots,k}}\right\|_2^2 = \left\|\sqrt{\varphi_{1,\ldots,n-1}}\right\|_2^2 + \left\|\sqrt{\varphi_{1,\ldots,k}}\right\|_2^2 - 2\int \sqrt{\varphi_{1,\ldots,n-1}}\sqrt{\varphi_{n,\ldots,k}}\,d\mathbb{P}$$
$$= 2\left(1 - \int \sqrt{\varphi_{1,\ldots,n-1}}\sqrt{\varphi_{n,\ldots,k}}\,d\mathbb{P}\right) \tag{6}$$

Now $\int \sqrt{\varphi_{1,\ldots,n-1}}\sqrt{\varphi_{n,\ldots,k}}\,d\mathbb{P} = \int \sqrt{\varphi_{1,\ldots,n-1}}\mathbb{E}\!\left(\sqrt{\varphi_{n,\ldots,k}} \mid \mathfrak{F}_n\right)d\mathbb{P} = \int \sqrt{\varphi_{1,\ldots,n-1}}M_{n,k}\,d\mathbb{P}$. Since $M_{n,k}$ is $\mathfrak{F}_n$-measurable, we have $\int \sqrt{\varphi_{1,\ldots,n-1}}M_{n,k}\,d\mathbb{P} = \int M_{n,k}\,d\mathbb{P}'$. Now since $M_{n,k}$ is increasing in $k$, we have that $2\left(1 - \int M_{n,k}\,d\mathbb{P}\right)$ is decreasing. Furthermore $\lim_n \lim_k 2\left(1 - \int M_{n,k}\,d\mathbb{P}'\right) = 0$ by dominated convergence. By Lemma 4.4.1, $\sqrt{\varphi_{1,\ldots,n}}$ is a Cauchy sequence and has a limit. $\square$

**Theorem 4.4.3**: Let $M = 1$ $(\mathbb{P} + \mathbb{P}')$-a.e. Then $\mathbb{P} \sim \mathbb{P}'$.

*Proof*: By symmetry it suffices to show that $\mathbb{P}' \ll \mathbb{P}$. Define $\varphi$ through $\sqrt{\varphi} = \lim_n \sqrt{\varphi_{1,\ldots,n}}$. By Lemma 4.4.2, this is well defined and $\sqrt{\varphi_{1,\ldots,n}} \to \sqrt{\varphi}$ in $L^2(\mathbb{P})$. By Hölder's inequality, $\sqrt{\varphi_{1,\ldots,n}}^2 \to \sqrt{\varphi}^2$ in $L^1(\mathbb{P})$. Therefore, $\varphi_{1,\ldots,n} \to \varphi$ in $L^1(\mathbb{P})$. We want to see that $\mathbb{P}' = \varphi \cdot \mathbb{P}$. Let $n \in \mathbb{N}$ and $A \in \mathfrak{F}_n$. Then for all $k > n$, we have $\mathbb{P}'(A) = \int_A \varphi_{1,\ldots,k}\,d\mathbb{P} \to \int_A \varphi\,d\mathbb{P}$. Therefore $\mathbb{P}'(A) = \int_A \varphi\,d\mathbb{P}$ for all $n \in \mathbb{N}$ and $A \in \mathfrak{F}_n$. The sets $A \in \mathfrak{F}$ that fulfill this equality form a Dynkin system. Since the equality holds all $A$ in the $\cap$-stable generator $\bigcup_{n\in\mathbb{N}}\mathfrak{F}_n$, it holds for all $A \in \mathfrak{F}$. $\square$

# 5 Conclusion

We state the theorem that we have proven in Section 4.

**Theorem 5.1**: Let $\mathbb{P}$ and $\mathbb{P}'$ be probability distributions on a filtered space $(\Omega, \mathfrak{F}, \mathfrak{F}_{n\in\mathbb{N}})$, s.t. $\mathfrak{F}_1 = \emptyset$, $\mathfrak{F} = \sigma\big(\bigcup_{n\in\mathbb{N}}\mathfrak{F}_n\big)$, and $\mathbb{P}|_{\mathfrak{F}_n} \sim \mathbb{P}'|_{\mathfrak{F}_n}$ for all $n \in \mathbb{N}$. For $n \in \mathbb{N}$, let $\Phi_n = \frac{\mathrm{d}\mathbb{P}'|_{\mathfrak{F}_n}}{\mathrm{d}\mathbb{P}|_{\mathfrak{F}_n}}$, and let $\varphi_n = \frac{\Phi_{n+1}}{\Phi_n}$. Then $M = \lim_n \lim_k \mathbb{E}\big(\prod_{i=n}^{k} \sqrt{\varphi_i} \mid \mathfrak{F}_n\big)$ exists $(\mathbb{P}+\mathbb{P}')$-a.e. and $\mathbb{P} \sim \mathbb{P}'$ if and only if $M = 1$ $(\mathbb{P}+\mathbb{P}')$-a.e. In this case, $(\Phi_n)^{1/2} \to (\frac{\mathrm{d}\mathbb{P}'}{\mathrm{d}\mathbb{P}})^{1/2}$ in $L^2(\mathbb{P})$

*Proof*: Existence of $M$ is proven in Section 4.1. For the second part, '$\Rightarrow$' is Corollary 4.3.2 and '$\Leftarrow$' is Theorem 4.4.3. The convergence of $\sqrt{\Phi_n}$ follows from Lemma 4.4.2 and Theorem 4.4.3. $\square$

We now apply the result to product spaces. Since the laws of families of random variables and stochastic processes are defined on a product space, this immediately gives a characterization of mutual absolute continuity of laws of families of random variables and laws of stochastic processes.

**Corollary 5.2**: Let $\mathbb{P}$ and $\mathbb{P}'$ be probability distributions on a product space $(\Omega, \mathcal{A}) = \big(\bigtimes_{i\in I}\Omega_i, \bigotimes_{i\in I}\mathcal{A}_i\big)$, s.t. $\mathbb{P}_{\pi_J} \sim \mathbb{P}'_{\pi_J}$ for all $J \subseteq I$ finite. For $J, K \subseteq I$ finite, let $\Phi_J = \frac{\mathrm{d}\mathbb{P}'_{\pi_J}}{\mathrm{d}\mathbb{P}_{\pi_J}}$, and let $\varphi_{J,K} = \frac{\Phi_J}{\Phi_K}$. Then for all increasing sequences $J = (J_n)_{n\in\mathbb{N}}$ of finite subsets of $I$ with $J_1 = \emptyset$, the random variable $M_J = \lim_n \lim_k \mathbb{E}\big(\prod_{i=n}^{k} \sqrt{\varphi_{J_i}} \mid \pi_{J_n}\big)$ exists $(\mathbb{P}+\mathbb{P}')$-a.e. Furthermore, $\mathbb{P} \sim \mathbb{P}'$ if and only if $M_J = 1$ $(\mathbb{P}+\mathbb{P}')$-a.e. for all increasing sequences $J = (J_n)_{n\in\mathbb{N}}$ of finite subsets of $I$ with $J_1 = \emptyset$. In this case $(\Phi_{J_n})^{1/2} \to (\frac{\mathrm{d}\mathbb{P}'_{\pi_K}}{\mathrm{d}\mathbb{P}_{\pi_K}})^{1/2}$, where $J_n$ is any sequence like above and $K = \bigcup_{n\in\mathbb{N}} J_n$.

*Proof*: First we prove that $\mathbb{P} \sim \mathbb{P}' \Leftrightarrow \forall J \subseteq I$ countable: $\mathbb{P}_{\pi_J} \sim \mathbb{P}'_{\pi_J}$. '$\Rightarrow$' is trivial. '$\Leftarrow$': By symmetry, it suffices to show $\mathbb{P}' \ll \mathbb{P}$ Let $A \in \mathcal{A}$, s.t. $\mathbb{P}(A) = 0$. Then there exists a countable $J \subseteq I$ and $\tilde{A}$, s.t. $A = \pi_J^{-1}\big(\tilde{A}\big)$. Therefore, $0 = \mathbb{P}(A) = \mathbb{P}_{\pi_J}\big(\tilde{A}\big) \Rightarrow 0 = \mathbb{P}'_{\pi_J}\big(\tilde{A}\big) = \mathbb{P}'(A)$.

The corollary now follows immediately by applying Theorem 5.1 with the filtration $(\mathfrak{F}_n)_{n\in\mathbb{N}}$ defined by $\mathfrak{F}_n = \sigma\big(\pi_{J_n}\big)$ for all increasing sequences $J = (J_n)_{n\in\mathbb{N}}$ of finite subsets of $I$ with $J_1 = \emptyset$. $\square$